# ON SMALLEST TRIANGLES

## GEOFFREY GRIMMETT AND SVANTE JANSON

ABSTRACT. Pick $n$ points independently at random in $\mathbb{R}^2$, according to a prescribed probability measure $\mu$, and let $\Delta_1^n \leq \Delta_2^n \leq \ldots$ be the areas of the $\binom{n}{3}$ triangles thus formed, in non-decreasing order. If $\mu$ is absolutely continuous with respect to Lebesgue measure, then, under weak conditions, the set $\{n^3 \Delta_i^n : i \geq 1\}$ converges as $n \to \infty$ to a Poisson process with a constant intensity $\kappa(\mu)$. This result, and related conclusions, are proved using standard arguments of Poisson approximation, and may be extended to functionals more general than the area of a triangle. It is proved in addition that, if $\mu$ is the uniform probability measure on the region $S$, then $\kappa(\mu) \leq 2/|S|$, where $|S|$ denotes the area of $S$. Equality holds in that $\kappa(\mu) = 2/|S|$ if $S$ is convex, and essentially only then. This work generalizes and extends considerably the conclusions of a recent paper of Jiang, Li, and Vitányi.

## 1. The problem

Drop $n$ points independently and uniformly at random into the unit square $[0,1]^2$ of $\mathbb{R}^2$. Of the $\binom{n}{3}$ triangles thus formed, let $\Delta_1^n$ be the smallest area. It is proved in [10] that $\mathbb{E}(\Delta_1^n)$ lies between $c_1 n^{-3}$ and $c_2 n^{-3}$ for certain positive real $c_1$ and $c_2$. An association is proposed in [10] between this problem and the (much harder) Heilbronn problem of finding the largest value of $\Delta_1^n$ for deterministic arrangements of $n$ points; see [3] and the references therein. A novel aspect of [10] is the use of Kolmogorov complexity.

Our target in this paper is to present a considerable strengthening of the above result, and to extend it to more general measures than uniform measure on $[0,1]^2$. This will be done using the now standard approach of Poisson approximation and the so-called Chen–Stein method, see [2]. From this point of view, the problem is a fairly simple exercise in probability theory, and does not of itself merit the intervention of concepts imported from complexity theory.

In broad terms, our extensions are as follows. First, we shall prove for the above problem on the square that $n^3 \Delta_1^n$ has, asymptotically, the exponential distribution with parameter 2 (that is, with mean $\frac{1}{2}$). Furthermore, the sequence $\Delta_1^n \leq \Delta_2^n \leq \ldots$ of smallest areas of triangles, written in non-decreasing order, is such that $n^3 \{\Delta_1^n, \Delta_2^n, \ldots\}$ converges as $n \to \infty$ to a Poisson process with intensity 2. These weak-limit theorems are complemented by the convergence of expectations obtained via an exponential bound on the tails of the $\Delta_i^n$.

Very similar results are valid when the underlying point distribution $\mu$ is uniform on a region $S$ ($\subseteq \mathbb{R}^2$) having non-zero finite Lebesgue measure $|S|$. In this case, the constant 2 is replaced by a constant $\kappa(\mu)$ satisfying $\kappa(\mu) \leq 2/|S|$. Furthermore, $\kappa(\mu) = 2/|S|$ if $S$ is convex, and more fully $\kappa(\mu) = 2/|S|$ if and only if $S$ differs from a convex set by a null set. Strict inequality holds, for example, if $S$ is the closure of its interior and is non-convex.

Indeed, corresponding results hold for a family of $n$ points drawn from $\mathbb{R}^2$ according to a general probability measure $\mu$ having a density function $f$ satisfying the boundedness conditions (A1) and (A2) below. In this very general situation, one obtains a formula for the ensuing constant $\kappa(\mu)$ in terms of $f$. Our results in this general setting are stated in Section 2. In Section 3, we set the scene by adapting to the present problem some of the basic theorems of Poisson approximation (see [2]) and exponential bounds (see [9]).

We return in Section 4 to the case of the uniform distribution on the unit square, where we give explicit proofs for this case. Special prominence is afforded to this instance since it has featured in [10], and since its resolution contains most of the basic ingredients of the general case. It is in addition useful to see the required arguments displayed in their simplicity without the intervention of certain minor but significant complications which arise in the setting of general measures.

Proofs in the general setting may be found in Section 5, and some properties of the constant $\kappa(\mu)$ are established in Section 6.

Further extensions come readily to mind, but we have not worked through the details. For example, for random points in $\mathbb{R}^3$, one might study the smallest area of a triangle formed by three points, or the smallest volume of a tetrahedron formed by four points, and similarly in higher dimensions. One might also study problems arising with points in other (Riemannian) manifolds, for example the smallest area of a geodetic triangle formed by random points on a torus or a sphere. Yet another possibility would be to allow $\mu$ to be singular, perhaps uniform measure on an arc of $\mathbb{R}^2$. The basic ingredients of such extensions are the arguments used in Section 3.

The work reported here is quite closely related to that of [14], where Poisson limit theorems were established for certain $U$-statistics, examples of which include the number of short inter-point distances and the number of 'flat' triangles created by an independent family of points in $\mathbb{R}^2$ with a common density function.

Our triangle problem has a property which raises its visibility above that of some other applications of the approach, namely the following. We will see in the proofs of Sections 4 and 5 that, although the smallest triangle has small area (with order $n^{-3}$), its diameter does not tend (stochastically) to 0 as $n \to \infty$, but remains of order 1. This fact underlies the observation that, in the case of the uniform distribution $\mu$ on a region $S$, the shape of $S$ affects the value of the constant $\kappa(\mu)$.

## 2. Notation and results

Let $\mu$ be a probability measure on the Borel $\sigma$-field of $\mathbb{R}^2$, and let $\mathbf{X}, \mathbf{X}_1, \mathbf{X}_2, \ldots, \mathbf{X}_n$ be points chosen independently at random according to the measure $\mu$. The $n$ points $\mathbf{X}_1, \mathbf{X}_2, \ldots, \mathbf{X}_n$ form $\binom{n}{3}$ triangles, whose areas we write in non-decreasing order as $\Delta_1^n \leq \Delta_2^n \leq \ldots$. For $\alpha > 0$, we shall be interested in the number $T^n(\alpha)$ of such triangles with area not exceeding $\alpha n^{-3}$, which is to say that

$$T^n(\alpha) = |\{k : \Delta_k^n \leq \alpha n^{-3}\}|.$$

The target of this paper is to establish conditions on $\mu$ under which $T^n(\alpha)$ has an approximate Poisson distribution, and more generally that the set $n^3\{\Delta_1^n, \Delta_2^n, \ldots\}$ is approximately distributed as a Poisson process.

We assume throughout that $\mu$ is absolutely continuous with respect to Lebesgue measure, which is to say that there exists an integrable (density) function $f : \mathbb{R}^2 \to [0, \infty)$ such that $\mu(A) = \int_A f(\mathbf{x}) \, d\mathbf{x}$ for any Borel set $A$. There follow two conditions on $f$.

(A1) *The function $f$ is bounded, in that there exists $M < \infty$ such that $f(\mathbf{x}) \leq M$ for all $\mathbf{x} \in \mathbb{R}^2$.*

This condition will prevent the clustering of points. A further condition concerns certain marginals of $f$, and this is prefaced with some notation. For notational convenience, we identify $\mathbb{R}^2$ with the complex plane, and we let $L(r, \theta)$ denote the line $\{re^{i\theta} + ite^{i\theta} : -\infty < t < \infty\}$. Note that $L(r, \theta) = \{\mathbf{x} : \mathbf{x} \cdot (\cos\theta, \sin\theta) = r\}$, and that $L(-r, \theta) = L(r, \theta + \pi)$. The family of all lines is given by either $\{L(r, \theta) : 0 \leq \theta < 2\pi, \ r \geq 0\}$, where $L(0, \theta)$ is identified with $L(0, \theta + \pi)$, or as $\{L(r, \theta) : 0 \leq \theta < \pi, \ -\infty < r < \infty\}$. We write $f_L$ for the integral of $f$ along the line $L$, and denote $f_{L(r,\theta)}$ by $f_{r,\theta}$ where convenient. Note that the mapping $r \mapsto f_{r,\theta}$ is the density function of the projection $\tau_\theta(\mathbf{X}) = \mathbf{X} \cdot (\cos\theta, \sin\theta)$ of $\mathbf{X}$ onto the line through the origin in the direction $\theta$. We shall assume that all these marginal densities are uniformly bounded.

(A2) *There exists $N < \infty$ such that, for all $r, \theta$, $f_{r,\theta} \leq N$.*

In our first theorem, we identify a constant which will play an important role in what follows. We write $|\mathbf{x}|$ for the Euclidean norm of $\mathbf{x}$ $(\in \mathbb{R}^2)$.

**Theorem 2.1.** *Assume that* (A1) *and* (A2) *hold, and let* $\mathbf{X}, \mathbf{Y}$ *be chosen independently according to* $\mu$. *Let* $L(R, \Theta)$ *be the straight line passing through the points* $\mathbf{X}$ *and* $\mathbf{Y}$. *The constant*

$$(2.1) \qquad \kappa = \kappa(\mu) = \frac{2}{3}\mathbb{E}\left(\frac{f_{R,\Theta}}{|\mathbf{X} - \mathbf{Y}|}\right)$$

*exists and satisfies* $0 < \kappa < \infty$. *Furthermore,*

$$(2.2) \qquad \kappa = \frac{2}{3}\int_0^\pi \int_{-\infty}^\infty f_{r,\theta}^3 \, dr \, d\theta.$$

We are ready to state our main theorem. For accounts of vague and weak convergence, see [8], [11].

**Theorem 2.2.** *Assume that* (A1) *and* (A2) *hold.*

(a) *The set* $n^3\{\Delta_1^n, \Delta_2^n, \ldots\}$ *converges vaguely to a Poisson process with constant intensity* $\kappa = \kappa(\mu)$, *and the convergence is weak when restricted to any* *bounded interval of $\mathbb{R}$. In particular, the random variable $n^3 \Delta_1^n$ converges weakly, as $n \to \infty$, to the exponential distribution with parameter $\kappa$.*

(b) *If $\mathbb{E}(|\mathbf{X}|^\delta) < \infty$ for some $\delta > 0$, then all positive moments of $n^3 \Delta_1^n$ converge to the corresponding moments of the exponential distribution with parameter $\kappa$. In particular, $\mathbb{E}(\Delta_1^n) \sim 1/(\kappa n^3)$ as $n \to \infty$.*

(c) *There exist absolute positive constants $A = A(\mu)$, $B = B(\mu)$ such that the total-variation distance between $T^n(\alpha)$ and the Poisson distribution with parameter $\lambda = \lambda^n(\alpha) = \mathbb{E}(T^n(\alpha))$ satisfies*

$$d_{\mathrm{TV}}\big(T^n(\alpha), \mathrm{Po}(\lambda^n(\alpha))\big) \le A \frac{1 - e^{-\lambda}}{\lambda} \cdot \frac{\alpha^2 \log(n^3/\alpha)}{n},$$

*for $0 < \alpha/n^3 < B$ and $n \ge 3$, where, for $\alpha > 0$, $\lambda = \lambda^n(\alpha) \to \kappa\alpha$ as $n \to \infty$.*

The conclusion of part (b) is valid with $\Delta_1^n$ replaced by $\Delta_{i+1}^n - \Delta_i^n$ for any given $i$. We shall return briefly to this point after the statement of Theorem 3.4.

With $\mathbf{X}$ chosen according to $\mu$, let $V = V(\mu)$ denote the covariance matrix of the 2-vector $\mathbf{X}$, and write $|V|$ for the determinant of $V$. If $\mathbb{E}(|\mathbf{X}|^2) = \infty$, we write $|V| = \infty$. The constant $\kappa$ may be calculated numerically using the formula in (2.2), and there is a universal lower bound on the product $\kappa(\mu)|V(\mu)|^{1/2}$. More important, there is an explicit formula for $\kappa$ valid whenever $\mu$ is uniform on a convex region. Prior to stating formally this last fact, we make a definition. A measurable subset $S$ of $\mathbb{R}^2$ is called *essentially convex* if there exists a set $N$ ($\subseteq \mathbb{R}^2$) with Lebesgue measure 0 such that the symmetric difference $S \triangle N$ is convex. It is trivial that every convex set is essentially convex.

**Theorem 2.3.** *Assume that* (A1) *and* (A2) *hold.*

(a) *The product $\kappa(\mu)|V(\mu)|^{1/2}$ is invariant under non-singular affine mappings of $\mathbb{R}^2$, and satisfies*

$$(2.3) \qquad\qquad \kappa(\mu)|V(\mu)|^{1/2} \ge \frac{1}{6\pi} \quad \text{for all } \mu.$$

(b) *Let $S$ be a subset of $\mathbb{R}^2$ with non-zero finite Lebesgue measure $|S|$, and let $\mu$ be the uniform probability measure on $S$. Then*

$$(2.4) \qquad\qquad\qquad\qquad \kappa(\mu) \le \frac{2}{|S|}.$$

*Equality holds in (2.4) if and only if $S$ is essentially convex.*

Proofs of the above results may be found in Sections 5 and 6. Here is an example of Theorem 2.3(a) in action. Suppose that $\mu$ is the standard bivariate normal distribution on $\mathbb{R}^2$. Using the rotational invariance of $\mu$, we have that

$$f_{r,\theta} = \frac{1}{\sqrt{2\pi}} e^{-r^2/2}.$$

By Theorem 2.1,

$$\kappa = \frac{2}{3}\pi \int_{-\infty}^{\infty} \frac{1}{(2\pi)^{3/2}} e^{-3r^2/2}\, dr = \frac{1}{3\sqrt{3}}.$$

By Theorem 2.3(a), if $\mu$ is a non-degenerate bivariate normal distribution with covariance matrix $V$, then

$$\kappa(\mu) = \frac{1}{3\sqrt{3|V|}}.$$

It is easy to see that, in contrast to (2.3), the quantity $\kappa(\mu)|V(\mu)|^{1/2}$ is not bounded above, even for distributions $\mu$ with $\mathbb{E}(|\mathbf{X}|^2) < \infty$. A simple example is given by the mixture of two symmetric bivariate normal distributions with covariance matrices $I$ and $a^2I$ where $I$ is the identity matrix and $0 < a < \infty$. We omit the calculations.

Theorem 2.3(b) has a striking reformulation in the language of geometrical probability. Let $S$ be a region of $\mathbb{R}^2$ having Lebesgue measure satisfying $0 < |S| < \infty$, and let $\mathbf{X}$, $\mathbf{Y}$ be points chosen independently and uniformly from $S$ according to the uniform probability measure on $S$. Let $l(\mathbf{X}, \mathbf{Y})$ be the one-dimensional Lebesgue measure (length) of the intersection with $S$ of the doubly infinite straight line through $\mathbf{X}$ and $\mathbf{Y}$. Note that, in the earlier notation, $f_{R,\Theta} = l(\mathbf{X}, \mathbf{Y})/|S|$.

**Corollary 2.4.** *We have that*

$$\mathbb{E}\left(\frac{l(\mathbf{X}, \mathbf{Y})}{|\mathbf{X} - \mathbf{Y}|}\right) \leq 3$$

*for all $S$ with $0 < |S| < \infty$. Equality holds if and only if $S$ is essentially convex.*

In the case of convex $S$, this exact calculation amounts to a reformulation of a classical result of M. Crofton in integral geometry, namely the following. Let $l(r, \theta)$ be the length of the intersection of $L(r, \theta)$ with a given convex region $S$. Crofton proved that

$$\int_0^\pi \int_{-\infty}^\infty l(r, \theta)^3\, dr\, d\theta = 3|S|^2.$$

His claim may be found in [5], with ancillary references at [4], [13]. Some related methods and conclusions of Crofton are summarised in [7, Section 4.13].

We close this section with a remark concerning the integrals

$$I(p, S) = \int_0^\pi \int_{-\infty}^\infty l(r, \theta)^p\, dr\, d\theta, \quad \text{for } p > 0,$$

where $S$ is assumed convex and $l$ is given as above. It is easily seen that $I(1, S) = \pi|S|$, and Crofton's formula is that $I(3) = 3|S|^2$. One may ask whether or not $I(p, S)$ depends only on $|S|$ for any other value of $p$. The answer is negative, and indeed $I(p, S)$ is not even invariant under affine area-preserving maps of $S$. As an illustration, let $S_\epsilon$ be the rectangle $[0, \epsilon^{-1}] \times [0, \epsilon]$. By elementary estimates, as $\epsilon \to 0$,

$$I(p, S_\epsilon) \to \begin{cases} 0 & \text{for } p \in (1, 3), \\ \infty & \text{for } p \in (0, 1) \cup (3, \infty). \end{cases}$$

## 3. Generalities

In advance of proving the theorems of the last section, we collect together certain sufficient conditions for Poisson-limit results to hold. We begin with the 'one-dimensional' distribution of $T^n(\alpha)$, and follow Barbour and Eagleson [1] as interpreted in [2, Thm 2.N].

Let $\mathbf{U}, \mathbf{V}, \mathbf{X}, \mathbf{Y}, \mathbf{Z}$ be chosen independently from $\mathbb{R}^2$ according to $\mu$, and let $\Delta(\mathbf{X}, \mathbf{Y}, \mathbf{Z})$ be the area of the triangle with vertices $\mathbf{X}, \mathbf{Y}, \mathbf{Z}$. We fix $\beta > 0$, and write

$$
\begin{aligned}
\pi(\beta) &= \mathbb{P}\big(\Delta(\mathbf{X}, \mathbf{Y}, \mathbf{Z}) \le \beta\big), \\
(3.1) \qquad \pi_1(\beta) &= \mathbb{P}\big(\Delta(\mathbf{X}, \mathbf{Y}, \mathbf{Z}) \le \beta, \ \Delta(\mathbf{X}, \mathbf{U}, \mathbf{V}) \le \beta\big), \\
\pi_2(\beta) &= \mathbb{P}\big(\Delta(\mathbf{X}, \mathbf{Y}, \mathbf{Z}) \le \beta, \ \Delta(\mathbf{X}, \mathbf{Y}, \mathbf{V}) \le \beta\big).
\end{aligned}
$$

We let $\alpha > 0$ and write

$$
\lambda^n(\alpha) = \binom{n}{3} \pi(\alpha n^{-3}),
$$

the mean number of triangles with area not exceeding $\alpha n^{-3}$. The following inequality will be useful.

**Lemma 3.1.** *We have that* $\pi(\beta)^2 \le \pi_1(\beta) \le \pi_2(\beta)$.

*Proof.* Let $H$ be the perpendicular distance (counted non-negative) from $\mathbf{Z}$ to the doubly infinite straight line through $\mathbf{X}$ and $\mathbf{Y}$, and write $\Delta = \Delta(\mathbf{X}, \mathbf{Y}, \mathbf{Z})$. Since $\Delta = \frac{1}{2} H |\mathbf{X} - \mathbf{Y}|$, we have that

$$
(3.2) \qquad \pi(\beta) = \mathbb{E}\left(\mathbb{P}\left(H \le \frac{2\beta}{|\mathbf{X} - \mathbf{Y}|} \,\Big|\, \mathbf{X}\right)\right),
$$

and the Cauchy–Schwarz inequality yields

$$
\begin{aligned}
\pi(\beta)^2 &= \mathbb{E}\left(\mathbb{P}\left(H \le \frac{2\beta}{|\mathbf{X} - \mathbf{Y}|} \,\Big|\, \mathbf{X}\right)\right)^2 \\
&\le \mathbb{E}\left(\mathbb{P}\left(H \le \frac{2\beta}{|\mathbf{X} - \mathbf{Y}|} \,\Big|\, \mathbf{X}\right)^2\right) = \pi_1(\beta).
\end{aligned}
$$

Similarly,

$$
\begin{aligned}
\pi_1(\beta) &= \mathbb{E}\left(\mathbb{E}\left(\mathbb{P}\left(H \le \frac{2\beta}{|\mathbf{X} - \mathbf{Y}|} \,\Big|\, \mathbf{X}, \mathbf{Y}\right) \,\Big|\, \mathbf{X}\right)^2\right) \\
&\le \mathbb{E}\left(\mathbb{E}\left(\mathbb{P}\left(H \le \frac{2\beta}{|\mathbf{X} - \mathbf{Y}|} \,\Big|\, \mathbf{X}, \mathbf{Y}\right)^2 \,\Big|\, \mathbf{X}\right)\right) = \pi_2(\beta),
\end{aligned}
$$

as required. $\qquad\square$

For random variables $A$, $B$, we define the *total-variation distance*

$$
d_{\mathrm{TV}}(A, B) = \sup_E \big|\mathbb{P}(A \in E) - \mathbb{P}(B \in E)\big|
$$

where the supremum is taken all events $E$. We write $\mathrm{Po}(\gamma)$ for a random variable having the Poisson distribution with parameter $\gamma$.

**Theorem 3.2.** *For $\alpha > 0$ and $n \geq 3$,*

$$d_{\mathrm{TV}}\big(T^n(\alpha), \mathrm{Po}(\lambda^n(\alpha))\big) \leq \frac{1 - e^{-\lambda}}{2\lambda} n^5 \pi_2,$$

*where $\lambda = \lambda^n(\alpha)$ and $\pi_2 = \pi_2(\alpha n^{-3})$.*

*Proof.* The set of areas of the $\binom{n}{3}$ triangles formed by the points $\mathbf{X}_1, \mathbf{X}_2, \ldots, \mathbf{X}_n$ is a dissociated family of random variables. We now apply [2, Thm 2.N] with $\beta = \alpha n^{-3}$. The upper bound follows via an application of Lemma 3.1. $\qquad\square$

It is an immediate corollary that

$$\mathbb{P}(\Delta_1^n > \alpha n^{-3}) = \mathbb{P}(T^n(\alpha) = 0)$$

satisfies

(3.3) $$\big|\mathbb{P}(\Delta_1^n > \alpha n^{-3}) - e^{-\lambda^n(\alpha)}\big| \leq \frac{1 - e^{-\lambda}}{2\lambda} n^5 \pi_2,$$

whence $n^3 \Delta_1^n$ has, asymptotically as $n \to \infty$, the exponential distribution with parameter $\lim_{n \to \infty} \lambda^n(\alpha)/\alpha$, whenever: (a) this limit exists and is independent of $\alpha$, and (b) the right-hand side of (3.3) tends to zero for every $\alpha > 0$. We shall verify these conditions for suitable $\mu$ in the subsequent sections.

We turn now to the set $n^3 \Delta = n^3 \{\Delta_1^n, \Delta_2^n, \ldots\}$.

**Theorem 3.3.** *Assume that there exist constants $c, c' \in (0, \infty)$ such that, for all $\alpha > 0$:*

(a) *$\lambda^n(\alpha) \to c\alpha$ as $n \to \infty$,*
(b) *$\lambda^n(\alpha) \leq c'\alpha$ for all $n$, and*
(c) *$n^5 \pi_2(\alpha n^{-3}) \to 0$ as $n \to \infty$.*

*Then $n^3 \Delta^n$ converges vaguely to a Poisson progress $\xi$ with constant intensity $c$, and, for all $T > 0$, the set $n^3 \Delta^n \cap [0, T]$ converges weakly to the process $\xi$ restricted to $[0, T]$. In particular, $n^3 \Delta_1^n$ converges weakly to the exponential distribution with parameter $c$.*

*Proof.* We apply Theorem 3.2 and Lemma 2.2 of [8] (see also [11]), with $\mathcal{I}$ the semi-ring of half-open intervals of the form $(a, b]$, $a, b \geq 0$. Note that, in the notation of [8], if $U \in \mathcal{I}$ is such that $U \subseteq [0, T]$ for some $T$, then

$$\mathbb{P}(X_j \in U) \leq \mathbb{P}(X_j \leq T),$$

and similarly for $\mathbb{P}(X_j \in U, \ X_k \in U')$. Lemma 3.1 is used to bound $\pi_1$ by $\pi_2$. $\qquad\square$

In Sections 4 and 5, we shall verify conditions (a), (b), (c) of Theorem 3.3 under the conditions stated in Section 2.

The above theorems will be used to prove weak convergence, but they can imply only little about the convergence of expectations. We concentrate here on the random variable $n^3 \Delta_1^n$, noting that similar arguments are valid for other elements of the set $n^3 \Delta^n$. We shall require an estimate for the tail of $n^3 \Delta_1^n$, in combination with the results above, in order to study its expectations.

**Theorem 3.4.** *We have that*

$$\mathbb{P}(n^3 \Delta_1^n > \alpha) \leq \exp\left(-\lambda^n(\alpha)^2/M_n(\alpha)\right) \quad \text{for } \alpha > 0$$

*where*

$$M_n(\alpha) = \max\left\{4\lambda^n(\alpha), 6n^5\pi_2(\alpha n^{-3})\right\}.$$

*Proof.* This is an application of Lemma 3.1, and Theorem 4 of [9]. □

A corresponding exponential bound for the tail of any given $\Delta_i^n$, for $i \geq 1$, follows from Theorem 10 of [9], and may be used in a similar way to that of the proof of Theorem 2.2(b) in order to obtain the convergence of the moments of $\Delta_i^n$ as $n \to \infty$; see the remark after the statement of Theorem 2.2.

## 4. Uniform distribution on the unit square

This section contains an analysis of the special instance in which $\mu$ is the uniform measure on the unit square $S = [0,1]^2$. We give special emphasis to this case for two reasons. First, it is the case considered in [10], and secondly, it is valuable to see an example of the required estimates where other minor but significant matters do not arise.

We assume henceforth in this section that $\mu$ is the uniform measure on the unit square $S$, and we next identify the constant $c$ for an application of Theorem 3.3.

**Theorem 4.1.** *Let $\mathbf{X}, \mathbf{Y}$ be chosen independently according to $\mu$, and let $l(\mathbf{X}, \mathbf{Y})$ be the Lebesgue measure (length) of the intersection with $S$ of the doubly infinite straight line through $\mathbf{X}$ and $\mathbf{Y}$. The constant*

$$\kappa = \frac{2}{3}\mathbb{E}\left(\frac{l(\mathbf{X}, \mathbf{Y})}{|\mathbf{X} - \mathbf{Y}|}\right) \tag{4.1}$$

*exists and satisfies $0 < \kappa < \infty$.*

We shall see in the proof of Theorem 2.3 (see Section 6) that $\kappa = 2$ for this special instance of $\mu$.

*Proof.* The basic estimate necessary for the proof is the following. There exists an absolute constant $\nu$ ($\geq 1$) such that

$$\mathbb{P}(|\mathbf{X} - \mathbf{Y}| \leq \epsilon \,|\, \mathbf{X}) \leq \nu\epsilon^2 \quad \text{for all } \epsilon > 0 \text{ and all } \mathbf{X}. \tag{4.2}$$

For reasons which will emerge soon, we shall assume that

$$\nu \geq e. \tag{4.3}$$

It follows that

$$\mathbb{E}(|\mathbf{X} - \mathbf{Y}|^{-1}) = \int_0^\infty \mathbb{P}(|\mathbf{X} - \mathbf{Y}|^{-1} > u)\, du \leq \int_0^\infty \min\{1, \nu u^{-2}\}\, du = 2\sqrt{\nu}.$$

The diameter of $S$ is $\sqrt{2}$, and therefore,

$$\kappa \leq \frac{2\sqrt{2}}{3}\mathbb{E}\left(|\mathbf{X} - \mathbf{Y}|^{-1}\right) = \frac{4\sqrt{2\nu}}{3},$$

as required. □

We next verify the conditions of Theorem 3.3 with $c = \kappa$, obtaining thereby the Poisson-process limit for the set $n^3\Delta^n$, as $n \to \infty$. Conditions (a) and (b) are implied by the next lemma, proved below.

**Lemma 4.2.** *We have that $\lambda^n(\alpha) \leq \kappa\alpha$ for all $\alpha, n$, and $\lambda^n(\alpha) \to \kappa\alpha$ as $n \to \infty$, for all $\alpha > 0$.*

We shall require an upper bound for $\pi_2(\beta)$ valid for small $\beta$.

**Lemma 4.3.** *There exist absolute constants $A$, $B$, with $B > 0$, such that $\pi_2(\beta) \leq A\beta^2 \log \beta^{-1}$ for $\beta < B$.*

It is an immediate consequence of Lemma 4.3 that

$$(4.4) \qquad \pi_2(\alpha n^{-3}) \leq \frac{A\alpha^2}{n^6} \log(n^3/\alpha) = \mathrm{O}(n^{-6} \log n) \quad \text{for } \alpha > 0 \text{ and large } n,$$

thus verifying condition (c) of Theorem 3.3. It follows that the set $n^3 \Delta^n$ satisfies the conclusions of Theorem 3.3, namely the vague and weak convergence to a Poisson-process limit. In addition, Theorem 3.2 provides a bound tending to zero for the total-variation distance between $T^n(\alpha)$ and the appropriate Poisson distribution. An extra ingredient is needed in order to prove convergence of expectations, as stated next.

**Theorem 4.4.** *For $p = 1, 2, \ldots$, the random variables $\{|n^3\Delta_1^n|^p : n \geq 3\}$ are uniformly integrable, and therefore, as $n \to \infty$, $\mathbb{E}(|n^3\Delta_1^n|^p) \to p! \, c^{-p}$, the $p$th moment of the exponential distribution with parameter $c$.*

The proof is given at the end of this section.

*Proof of Lemma 4.2.* The letter $L$ denotes a doubly infinite straight line of $\mathbb{R}^2$. For $\mathbf{x}, \mathbf{y} \in \mathbb{R}^2$, $\mathbf{x} \neq \mathbf{y}$, we write $L(\mathbf{x}, \mathbf{y})$ for the doubly-infinite straight line through $\mathbf{x}$ and $\mathbf{y}$. For $\epsilon > 0$, we write $S_\epsilon(\mathbf{x}, \mathbf{y})$ for the closed strip containing all points whose Euclidean distance to $L(\mathbf{x}, \mathbf{y})$ does not exceed $\epsilon$.

Let $\mathbf{X}$, $\mathbf{Y}$, $\mathbf{Z}$ be chosen independently according to $\mu$, and let $H$ be the absolute value of the length of the perpendicular dropped from $\mathbf{Z}$ onto $L(\mathbf{X}, \mathbf{Y})$. We have that $\Delta = \Delta(\mathbf{X}, \mathbf{Y}, \mathbf{Z})$ satisfies

$$(4.5) \qquad \mathbb{P}(\Delta \leq \beta \mid \mathbf{X}, \mathbf{Y}) = \mathbb{P}\left( H \leq \frac{2\beta}{|\mathbf{X} - \mathbf{Y}|} \,\bigg|\, \mathbf{X}, \mathbf{Y} \right)$$

$$= \mathbb{P}\Big( \mathbf{Z} \in S_{2\beta/|\mathbf{X}-\mathbf{Y}|}(\mathbf{X}, \mathbf{Y}) \,\Big|\, \mathbf{X}, \mathbf{Y} \Big) \leq \frac{4\beta l(\mathbf{X}, \mathbf{Y})}{|\mathbf{X} - \mathbf{Y}|}$$

with equality if the event $C(\beta)$ does not occur, where $C(\beta)$ is the event that $S_{2\beta/|\mathbf{X}-\mathbf{Y}|}(\mathbf{X}, \mathbf{Y})$ contains one or more of the four corners of $S$, and $1_E$ denotes the indicator function of an event $E$. Thus,

$$\frac{1}{\beta} \left| \mathbb{P}(\Delta \leq \beta) - \mathbb{E}\left( \frac{4\beta l(\mathbf{X}, \mathbf{Y})}{|\mathbf{X} - \mathbf{Y}|} \right) \right| \leq \mathbb{E}\left( \frac{4l(\mathbf{X}, \mathbf{Y})}{|\mathbf{X} - \mathbf{Y}|} 1_{C(\beta)} \right).$$

We have by the absolute continuity of $\mu$ that $\mathbb{P}(C(\beta)) \to 0$ as $\beta \to 0$, and it follows by Theorem 4.1 and the monotone convergence theorem that

$$(4.6) \qquad \frac{1}{\beta}\mathbb{P}(\Delta \leq \beta) \to \mathbb{E}\left( \frac{4l(\mathbf{X}, \mathbf{Y})}{|\mathbf{X} - \mathbf{Y}|} \right) \quad \text{as } \beta \downarrow 0.$$

As an immediate consequence, we obtain that

$$\lambda^n(\alpha) = \binom{n}{3} \mathbb{P}(\Delta \leq \alpha n^{-3}) = (1 + \mathrm{o}(1)) \frac{n^3}{6} \cdot \frac{\alpha}{n^3} \mathbb{E}\left(\frac{4l(\mathbf{X}, \mathbf{Y})}{|\mathbf{X} - \mathbf{Y}|}\right)$$
$$\to \kappa\alpha \quad \text{as } n \to \infty.$$

The uniform upper bound on $\lambda^n(\alpha)$ follows from (4.5). $\qquad \square$

*Proof of Lemma 4.3.* We have by (4.2) and (4.5) that, for $\epsilon > 0$,

$$(4.7) \qquad \pi_2(\beta) = \mathbb{E}\left(\mathbb{P}(\Delta \leq \beta \mid \mathbf{X}, \mathbf{Y})^2\right)$$
$$\leq \mathbb{E}\left(1_{\{|\mathbf{X}-\mathbf{Y}|\leq\epsilon\}} + \left(\frac{4\beta l(\mathbf{X}, \mathbf{Y})}{|\mathbf{X} - \mathbf{Y}|}\right)^2 1_{\{|\mathbf{X}-\mathbf{Y}|>\epsilon\}}\right)$$
$$\leq \nu\epsilon^2 + 32\beta^2 \mathbb{E}\left(\frac{1_{\{|\mathbf{X}-\mathbf{Y}|>\epsilon\}}}{|\mathbf{X} - \mathbf{Y}|^2}\right).$$

Now, by (4.2) and (4.3),

$$(4.8) \qquad \mathbb{E}\left(\frac{1_{\{|\mathbf{X}-\mathbf{Y}|>\epsilon\}}}{|\mathbf{X} - \mathbf{Y}|^2}\right) = \int_0^{\epsilon^{-2}} \mathbb{P}\left(|\mathbf{X} - \mathbf{Y}|^{-2} > u\right) du$$
$$\leq \int_0^{\epsilon^{-2}} \min\{1, \nu u^{-1}\} du$$
$$= \nu + \nu \log(\epsilon^{-2}/\nu) \leq \nu \log(\epsilon^{-2}),$$

if $\nu\epsilon^2 \leq 1$. We substitute this into (4.7) and minimize over $\epsilon$ to obtain that there exist $A = A(\nu)$ and $B = B(\nu) > 0$ such that

$$\pi_2(\beta) \leq A\beta^2 \log \beta^{-1} \quad \text{for} \quad \beta < B. \qquad \square$$

*Proof of Theorem 4.4.* By a standard argument (see [7, Section 7.10]), it suffices to show that, for each $p \geq 1$, the set $\{\mathbb{E}(|n^3 \Delta_1^n|^p) : n \geq 3\}$ is uniformly bounded. This will be achieved by the exponential tail estimate of Theorem 3.4, which we make more precise as follows.

Consider the expression $M_n(\alpha)$ in Theorem 3.4. By Lemma 4.2 and equation (4.4),

$$M_n(\alpha) \leq \max\left\{c_1\alpha, \frac{c_2\alpha^2 \log n}{n}\right\},$$

for positive constants $c_i$ and $n \geq 3$. Therefore,

$$M_n(\alpha) \leq \begin{cases} c_1\alpha & \text{if } \alpha < \dfrac{c_1 n}{c_2 \log n}, \\ \dfrac{c_2\alpha^2 \log n}{n} & \text{otherwise.} \end{cases}$$

Furthermore, by (4.6), $\mathbb{P}(\Delta \leq \beta) \geq c_3\beta$ for $0 < \beta < 1$, and thus $\lambda_n(\alpha) > c_4\alpha$ for $\alpha < n^3$ and $n \geq 3$. We apply Theorem 3.4 to obtain

$$(4.9) \qquad \mathbb{P}(n^3 \Delta_1^n > \alpha) \leq \begin{cases} e^{-c_5\alpha} & \text{if } \alpha < \dfrac{c_1 n}{c_2 \log n}, \\[2ex] \exp\left(-\dfrac{c_6 n}{\log n}\right) & \text{otherwise;} \end{cases}$$

this probability is evidently 0 for $\alpha \geq n^3$. Therefore,

$$(4.10) \qquad \mathbb{E}\left(|n^3 \Delta_1^n|^p\right) \leq \int_0^{c_7 n/\log n} p\alpha^{p-1} e^{-c_5\alpha}\, d\alpha$$
$$+ \int_{c_7 n/\log n}^{n^3} p\alpha^{p-1} \exp\left(-\frac{c_6 n}{\log n}\right) d\alpha,$$

which is bounded above uniformly in $n$. □

## 5. Proofs of Theorems 2.1 and 2.2

We assume throughout this section that (A1) and (A2) hold, and begin with a preliminary lemma. As in Section 2, a line in $\mathbb{R}^2$ is denoted $L(r, \theta)$, and $L(R, \Theta)$ denotes the line passing through two points $\mathbf{X}$, $\mathbf{Y}$ chosen independently according to the measure $\mu$.

**Lemma 5.1.** *Let $g(r, \theta)$ be a non-negative measurable function on the set of lines $L(r, \theta)$ in $\mathbb{R}^2$. Then*

$$\mathbb{E}\left(\frac{g(R, \Theta)}{|\mathbf{X} - \mathbf{Y}|}\right) = \int_0^\pi \int_{-\infty}^\infty g(r, \theta) f_{r,\theta}^2\, dr\, d\theta.$$

*Proof.* Using polar coordinates in the form $\mathbf{Y} = \mathbf{X} + ire^{i\theta}$, $r \in (-\infty, \infty)$, $\theta \in [0, \pi)$, we find

$$\mathbb{E}\left(\frac{g(R, \Theta)}{|\mathbf{X} - \mathbf{Y}|}\,\Big|\, \mathbf{X}\right) = \int_0^\pi \int_{-\infty}^\infty \frac{g(\tau_\theta(\mathbf{X}), \theta)}{|r|} f(\mathbf{X} + ire^{i\theta})|r|\, dr\, d\theta$$
$$= \int_0^\pi g(\tau_\theta(\mathbf{X}), \theta) f_{\tau_\theta(\mathbf{X}), \theta}\, d\theta,$$

where $\tau_\theta(\mathbf{X}) = \mathbf{X} \cdot (\cos\theta, \sin\theta)$ as before. Hence, by Fubini's theorem,

$$\mathbb{E}\left(\frac{g(R, \Theta)}{|\mathbf{X} - \mathbf{Y}|}\right) = \int_0^\pi \mathbb{E}\left(g(\tau_\theta(\mathbf{X}), \theta) f_{\tau_\theta(\mathbf{X}), \theta}\right)\, d\theta = \int_0^\pi \int_{-\infty}^\infty g(r, \theta) f_{r,\theta}^2\, dr\, d\theta,$$

since $\tau_\theta(\mathbf{X})$ has density function $f_{r,\theta}$. □

*Proof of Theorem 2.1.* It is elementary that $\kappa$ exists and is non-negative. Furthermore,

$$\kappa \leq \frac{2N}{3}\mathbb{E}\big(|\mathbf{X} - \mathbf{Y}|^{-1}\big).$$

The finiteness of $\mathbb{E}\big(|\mathbf{X} - \mathbf{Y}|^{-1}\big)$ follows just as in the proof of Theorem 4.1, or alternatively as follows. By Lemma 5.1,

$$(5.1) \qquad \mathbb{E}\big(|\mathbf{X} - \mathbf{Y}|^{-1}\big) = \int_0^\pi \int_{-\infty}^\infty f_{r,\theta}^2 \, dr \, d\theta \leq N \int_0^\pi \int_{-\infty}^\infty f_{r,\theta} \, dr \, d\theta = N\pi.$$

Equation (2.2) follows from Lemma 5.1 with $g(r,\theta) = f_{r,\theta}$.     $\square$

*Proof of Theorem 2.2.* We show first that conditions (a), (b), (c) of Theorem 3.3 are valid with $c = \kappa$. Let $\Delta = \Delta(\mathbf{X}, \mathbf{Y}, \mathbf{Z})$ where $\mathbf{X}$, $\mathbf{Y}$, $\mathbf{Z}$ are chosen independently according to $\mu$, and let $H$ be the length of the perpendicular dropped from $\mathbf{Z}$ onto $L(R, \Theta)$. Then

$$(5.2) \qquad \mathbb{P}(H \leq h \mid \mathbf{X}, \mathbf{Y}) = \int_{R-h}^{R+h} f_{r,\Theta} \, dr \leq 2hN.$$

Consequently, by (3.2) and (5.1),

$$\pi(\beta) \leq 4\beta N\mathbb{E}\big(|\mathbf{X} - \mathbf{Y}|^{-1}\big) \leq 4\beta N^2\pi,$$

and (b) follows for suitable $c'$.

For each $\theta$, the mapping $r \mapsto f_{r,\theta}$ is integrable. Let $E_\theta$ be the set of $r$ ($\in \mathbb{R}$) such that

$$\lim_{\epsilon \to 0} \frac{1}{2\epsilon} \int_{r-\epsilon}^{r+\epsilon} f_{s,\theta} \, ds = f_{r,\theta},$$

and note that, by the Lebesgue differentiation theorem (see [12], Theorem 9.3), the complement $\mathbb{R} \setminus E_\theta$ has Lebesgue measure zero. We thus have

$$(5.3) \qquad \frac{1}{\beta}\mathbb{P}\left(H \leq \frac{2\beta}{|\mathbf{X} - \mathbf{Y}|} \,\Big|\, \mathbf{X}, \mathbf{Y}\right) = \frac{1}{\beta} \int_{R-2\beta/|\mathbf{X}-\mathbf{Y}|}^{R+2\beta/|\mathbf{X}-\mathbf{Y}|} f_{r,\theta} \, dr$$
$$\to \frac{4f_{R,\Theta}}{|\mathbf{X} - \mathbf{Y}|} \quad \text{as } \beta \to 0,$$

whenever $R \notin E_\Theta$.

For $\mathbf{x}, \mathbf{y} \in \mathbb{R}^2$, $\mathbf{x} \neq \mathbf{y}$, write $L(r, \theta)$ for the line through $\mathbf{x}$ and $\mathbf{y}$. Define $u$, $v$ by $\mathbf{x} = (r + iu)e^{i\theta}$, $\mathbf{y} = (r + iv)e^{i\theta}$. The mapping $(\mathbf{x}, \mathbf{y}) \mapsto (r, \theta, u, v)$, with domain $\{(\mathbf{x}, \mathbf{y}) : \mathbf{x} \neq \mathbf{y}\}$ and range $\{(r, \theta, u, v) : u \neq v\}$, is differentiable with differentiable inverse. Thus it maps null sets to null sets. It follows that (5.3) holds for almost every $(\mathbf{X}, \mathbf{Y}) \in \mathbb{R}^2$, and thus $\mu$-a.s. Together with (5.2) and the dominated convergence theorem, we deduce that

$$\frac{1}{\beta}\pi(\beta) = \frac{1}{\beta}\mathbb{E}\left(\mathbb{P}\left(H \leq \frac{2\beta}{|\mathbf{X} - \mathbf{Y}|} \,\Big|\, \mathbf{X}, \mathbf{Y}\right)\right) \to 4\mathbb{E}\left(\frac{f_{R,\Theta}}{|\mathbf{X} - \mathbf{Y}|}\right) \quad \text{as } \beta \downarrow 0.$$

Consequently, for any $\alpha > 0$,

$$\lambda^n(\alpha) = \binom{n}{3} \pi(\alpha n^{-3}) \to \frac{4}{6} \mathbb{E}\left(\frac{f_{R,\Theta}}{|\mathbf{X} - \mathbf{Y}|}\right) \alpha = \kappa\alpha \quad \text{as } n \to \infty,$$

thus confirming (a).

Finally, (c) follows as in the proof of (4.4) and Lemma 4.3. Theorem 2.2(a) follows.

We turn now to Theorem 2.2(b). The proof of Theorem 4.4 may be followed until equation (4.9), which is valid subject to an upper bound on $\alpha/n^3$ (with $c_6$ depending on the bound). Unlike the situation when $\mu$ has bounded support, there is in the general case no deterministic upper bound on the maximal area of a triangle. Thus we shall need a further estimate on the tail of $\Delta_1^n$. Assume that $\delta \in (0, 1)$ is such that $\mathbb{E}(|\mathbf{X}|^\delta) < \infty$, and note that $\Delta = \Delta(\mathbf{X}, \mathbf{Y}, \mathbf{Z})$ satisfies

$$\Delta \leq |\mathbf{X} - \mathbf{Y}| \cdot |\mathbf{X} - \mathbf{Z}|,$$

whence

$$\mathbb{E}(\Delta^{\delta/2}) \leq \sqrt{\mathbb{E}(|\mathbf{X} - \mathbf{Y}|^\delta \cdot \mathbb{E}(|\mathbf{X} - \mathbf{Z}|^\delta)} = \mathbb{E}(|\mathbf{X} - \mathbf{Y}|^\delta) \leq 2\mathbb{E}(|\mathbf{X}|^\delta) < \infty.$$

By Markov's inequality, there exists $t_0 \geq 1$ such that

$$\mathbb{P}(\Delta > t) \leq t^{-\delta/2}\mathbb{E}(\Delta^{\delta/2}) \leq \tfrac{1}{8}t^{-\delta/3} \quad \text{for } t \geq t_0.$$

Therefore, by considering the $\lfloor n/3 \rfloor$ triangles $\mathbf{X}_1\mathbf{X}_2\mathbf{X}_3$, $\mathbf{X}_4\mathbf{X}_5\mathbf{X}_6$, ... with disjoint vertex sets,

$$\mathbb{P}(\Delta_1^n > t) \leq (8t^{\delta/3})^{-\lfloor n/3 \rfloor} \leq 2^{2-n}t^{-\lfloor n/3 \rfloor \delta/3}, \quad t \geq t_0.$$

Let $p \geq 1$. It follows from the above that

$$\mathbb{P}(\Delta_1^n > t) \leq 2^{2-n}t^{-p-1} \quad \text{for } n \geq 3 + \frac{9(p+1)}{\delta}, \ t \geq t_0.$$

Hence, for $n \geq 3 + 9(p+1)/\delta$,

$$\int_{t_0 n^3}^\infty p\alpha^{p-1}\mathbb{P}(n^3\Delta_1^n > \alpha) \, d\alpha = n^{3p} \int_{t_0}^\infty pt^{p-1}\mathbb{P}(\Delta_1^n > t) \, dt \leq n^{3p} \cdot 2^{2-n}p.$$

We combine this with (4.10), with the range of the second integral of that equation extended to $t_0 n^3$, to deduce as required that the $p$th moment of $\Delta_1^n$ is uniformly bounded in $n$. This completes the proof of Theorem 2.2(b).

Part (c) follows by Theorem 3.2 and the proof of Lemma 4.3. □

We remark that conditions (A1) and (A2) may be relaxed, but we do not attempt here to obtain the weakest necessary assumptions. We note as follows, however, that conditions (a) and (b) of Theorem 3.3, with $c = \kappa$, hold whenever $\kappa < \infty$. Under the assumption $\kappa < \infty$, we have by (2.2) and [15, Thm 1], that the Hardy–Littlewood maximal function of $f_{r,\theta}$ is bounded in $L^3$. It is now a simple matter to replace (5.2) in the above proof.

## 6. Proof of Theorem 2.3

Finally, we prove Theorem 2.3, beginning with part (a). Let $A$ be a non-singular $2 \times 2$ matrix, let $\mathbf{b} \in \mathbb{R}^2$, and write $\widetilde{\mathbf{X}}_j = A\mathbf{X}_j + \mathbf{b}$. Denote by $\widetilde{\mu}$ the measure governing the $\widetilde{\mathbf{X}}_j$. Let the areas formed by triangles of the $\widetilde{\mathbf{X}}_j$ be written $\widetilde{\Delta}_i^n$, so that $\widetilde{\Delta}_i^n = |A||\Delta_i^n|$. We have by Theorem 2.2 that $n^3\widetilde{\Delta}^n$ converges to a Poisson process with intensity $|A|^{-1}\kappa(\mu)$, and thus $\kappa(\widetilde{\mu}) = |A|^{-1}\kappa(\mu)$. The covariance matrix of $\widetilde{\mathbf{X}}_1$ has determinant $|V(\widetilde{\mu})| = |A|^2 \cdot |V(\mu)|$, and thus $\kappa(\widetilde{\mu})|V(\widetilde{\mu})|^{1/2} = \kappa(\mu)|V(\mu)|^{1/2}$ as required.

We prove (2.3) next. Since the mapping $r \mapsto f_{r,\theta}$ is a density function,

$$\int_0^\pi \int_{-\infty}^\infty f_{r,\theta}\, dr\, d\theta = \pi.$$

The Cauchy–Schwarz inequality yields

$$\left(\int_0^\pi \int_{-\infty}^\infty f_{r,\theta}^2\, dr\, d\theta\right)^2 \le \int_0^\pi \int_{-\infty}^\infty f_{r,\theta}\, dr\, d\theta \int_0^\pi \int_{-\infty}^\infty f_{r,\theta}^3\, dr\, d\theta = \pi \cdot \frac{3\kappa}{2},$$

by Theorem 2.1, whence, by (5.1),

(6.1)  $$\kappa(\mu) \ge \frac{2}{3\pi}\mathbb{E}\big(|\mathbf{X} - \mathbf{Y}|^{-1}\big)^2.$$

On the other hand, by Hölder's inequality with exponents $\frac{3}{2}$ and 3,

$$1 = \mathbb{E}\left(\frac{1}{|\mathbf{X} - \mathbf{Y}|^{2/3}} \cdot |\mathbf{X} - \mathbf{Y}|^{2/3}\right) \le \mathbb{E}\left(\frac{1}{|\mathbf{X} - \mathbf{Y}|}\right)^{2/3}\mathbb{E}\big(|\mathbf{X} - \mathbf{Y}|^2\big)^{1/3},$$

and thus

(6.2)  $$\mathbb{E}\big(|\mathbf{X} - \mathbf{Y}|^{-1}\big)^{-2} \le \mathbb{E}\big(|\mathbf{X} - \mathbf{Y}|^2\big) = 2(v_{11} + v_{22}),$$

where $V(\mu) = (v_{ij})$. If $V(\mu) = I$, the identity matrix, (6.1) and (6.2) yield

$$\kappa(\mu) \ge \frac{2}{3\pi} \cdot \frac{1}{2(v_{11} + v_{22})} = \frac{1}{6\pi} = \frac{1}{6\pi|V(\mu)|^{1/2}},$$

and (2.3) is proved in this case. The general case follows since $\kappa(\mu)|V(\mu)|^{1/2}$ is invariant under non-singular affine maps of $\mathbb{R}^2$.

We turn to part (b). First suppose that $S$ is convex. For $\mathbf{x} \in S$ and $\theta \in [0, 2\pi)$, let $r(\mathbf{x}, \theta) = \sup\{t : \mathbf{x} + te^{i\theta} \in S\}$. Then, for every $\mathbf{x} \in S$,

(6.3)  $$|S| = \int_0^{2\pi} \int_0^{r(\mathbf{x},\theta)} \rho\, d\rho\, d\theta = \int_0^{2\pi} \tfrac{1}{2}r(\mathbf{x}, \theta)^2\, d\theta,$$

and thus, by Fubini's theorem,

$$|S|^2 = \int_0^{2\pi} \int_S \tfrac{1}{2}r(\mathbf{x}, \theta)^2\, d\mathbf{x}\, d\theta.$$

The line $L(s, \theta)$ intersects $S$ in a (possibly empty) interval of length $l(s, \theta)$, say, and, for $\mathbf{x}$ in this interval, $r(\mathbf{x}, \theta + \frac{1}{2}\pi)$ varies linearly between 0 and $l(s, \theta)$. Hence,

$$\int_S r(\mathbf{x}, \theta + \tfrac{1}{2}\pi)^2 \, d\mathbf{x} = \int_{-\infty}^{\infty} \int_0^{l(s,\theta)} r^2 \, dr \, ds = \int_{-\infty}^{\infty} \tfrac{1}{3} l(s, \theta)^3 \, ds.$$

Consequently,

$$|S|^2 = \tfrac{1}{6} \int_0^{2\pi} \int_{-\infty}^{\infty} l(s, \theta)^3 \, ds \, d\theta = \tfrac{1}{3} \int_0^{\pi} \int_{-\infty}^{\infty} l(s, \theta)^3 \, ds \, d\theta,$$

and thus, since $f_{s,\theta} = l(s, \theta)/|S|$,

$$\kappa = \frac{2}{3} \int_0^{\pi} \int_{-\infty}^{\infty} f_{s,\theta}^3 \, ds \, d\theta = \frac{2}{|S|}.$$

The same conclusion is easily seen to hold under the weaker assumption that $S$ be *essentially* convex.

For a general region $S$, let $A(\mathbf{x}, \theta) = \{t \geq 0 : \mathbf{x} + te^{i\theta} \in S\}$. Then, as in (6.3),

$$(6.4) \qquad\qquad |S| \geq \int_0^{2\pi} \tfrac{1}{2} |A(\mathbf{x}, \theta)|^2 \, d\theta,$$

and, as above (except that $L(s, \theta) \cap S$ need not be an interval),

$$\int_0^{2\pi} \int_S \tfrac{1}{2} |A(\mathbf{x}, \theta)|^2 \, d\mathbf{x} \, d\theta = \tfrac{1}{3} \int_0^{\pi} \int_{-\infty}^{\infty} |L(s, \theta) \cap S|^3 \, ds \, d\theta$$

$$= \tfrac{1}{3} |S|^3 \int_0^{\pi} \int_{-\infty}^{\infty} f_{s,\theta}^3 \, ds \, d\theta,$$

yielding as required that $|S|^2 \geq \tfrac{1}{2} |S|^3 \kappa$. Equality holds only if there is equality in (6.4) for almost every $\mathbf{x}$, that is, if $A(\mathbf{x}, \theta)$, for almost every $(\mathbf{x}, \theta)$, differs only by a null set from an interval containing 0.

Let $T'$ denote the set of density points of a measurable subset $T$ of $\mathbb{R}^2$ (recall that $\mathbf{x}$ is a *density point* of $T$ if $|T \cap B(\mathbf{x}, r)|/|B(\mathbf{x}, r)| \to 1$ as $r \to 0$, where $B(\mathbf{x}, r)$ denotes the ball of radius $r$ with centre $\mathbf{x}$). It may be deduced from the above statement concerning the $A(\mathbf{x}, \theta)$ that $S'$ is a convex set (this is left to the reader to prove). The symmetric difference $T \triangle T'$ is a null set for any measurable $T$ (see [15, Section 1.1] and [6, Thm 2.9.11]), and in particular $N = S \triangle S'$ is a null set. Thus $S \triangle N$ is convex where $N$ is null. Hence $S$ is essentially convex and the proof is complete.

## Acknowledgements

GRG thanks John Marstrand for drawing his attention to this problem. Some of this work was done at EURANDOM, during the attendance by the authors at a meeting on discrete probability. David Stirzaker kindly indicated the relevance of the work of Morgan Crofton.

STATISTICAL LABORATORY, UNIVERSITY OF CAMBRIDGE, WILBERFORCE ROAD, CAMBRIDGE CB3 0WB, UNITED KINGDOM
   *E-mail*: g.r.grimmett@statslab.cam.ac.uk
   *URL*: http://www.statslab.cam.ac.uk/∼grg/

DEPARTMENT OF MATHEMATICS, UPPSALA UNIVERSITY, BOX 480, 751 06 UPPSALA, SWEDEN
   *E-mail*: svante@math.uu.se
   *URL*: http://www.math.uu.se/∼svante/